\documentclass[11pt, reqno]{amsart}
\usepackage{amsmath}
\usepackage{amssymb}
\input amssym.def
\input amssym.tex

\font\tengothic=eufm10
\font\sevengothic=eufm7
\newfam\gothicfam
      \textfont\gothicfam=\tengothic
      \scriptfont\gothicfam=\sevengothic
\def\goth#1{{\fam\gothicfam #1}}

\numberwithin{equation}{section}
\newenvironment{mydef}{{\bf Definition.}}{\hspace*{\fill} \par\vspace{1ex}}

\begin{document}
\setlength{\baselineskip}{1.7em}
\newtheorem{thm}{\bf Theorem}[section]
\newtheorem{pro}[thm]{\bf Proposition}
\newtheorem{claim}[thm]{\bf Claim}
\newtheorem{lemma}[thm]{\bf Lemma}
\newtheorem{cor}{\bf Corollary}[thm]

\newcommand{\mP}{{\mathbb P}}
\newcommand{\z}{{\mathbb Z}}
\newcommand{\nn}{{\mathbb N}}
\newcommand{\R}{{\mathcal R}}
\newcommand{\A}{{\mathcal A}}
\newcommand{\x}{{\mathbb X}}
\newcommand{\y}{{\mathbb Y}}
\newcommand{\ix}{I_{\mathbb X}}
\newcommand{\mi}{{\mathcal I}}
\newcommand{\B}{{\bf B}}
\newcommand{\bi}{{\bf I}}
\newcommand{\V}{{\bf V}}
\newcommand{\cv}{{\mathcal V}}
\newcommand{\bL}{{\bf L}}
\newcommand{\cl}{{\mathcal L}}
\newcommand{\cm}{{\mathcal M}}
\newcommand{\cn}{{\mathcal N}}
\newcommand{\m}{-\!\!--\!\!\rightarrow}
\newcommand{\smap}{\rightarrow\!\!\!\!\!\rightarrow}
\newcommand{\sfrac}[2]{\frac{\displaystyle #1}{\displaystyle #2}}
\newcommand{\gr}{\Gamma_{d+1}}
\newcommand{\img}{\Lambda_{d+1}}
\newcommand{\under}[1]{\underline{#1}}
\newcommand{\ov}[1]{\overline{#1}}
\newcommand{\lex}{{\le}_{\mbox{\scriptsize lex}}}

\title{Box-shaped matrices and the defining ideal of certain blowup surfaces}
\author{H\`a Huy T\`ai}
\thanks{Current address: Institute of Mathematics, P.O. Box 631 B\`o H\^o, H\`a N\^oi 10000, Vietnam. Email: tai@hanimath.ac.vn}
\address{Department of Mathematics and Statistics \\ Queen's University, Kingston ON K7L 3N6, Canada} \email{haht@mast.queensu.ca}
\subjclass[2000]{Primary: 13C40, 14J26, Secondary: 14E25.}
\begin{abstract}
In this paper, we generalize the notions of a matrix and its ideal of $2\times2$ minors to that of a box-shaped matrix and its ideal of $2\times2$ minors, and make use of these notions to study projective embeddings of certain blowup surfaces. We prove that the ideal of $2\times2$ minors of a generic box-shaped matrix is a perfect prime ideal that gives the algebraic description for the Segre embedding of the product of several projective spaces. We use the notion of the ideal of $2\times2$ minors of a box-shaped matrix to give an explicit description for the defining ideal of the blowup of $\mP^2$ along a set of ${d+1 \choose 2} ~ (d \in \z)$ points in generic position, embedded into projective spaces using very ample divisors which correspond to the linear systems of plane curves going through these points. 
\end{abstract}
\maketitle

\setcounter{section}{-1}
\section{Introduction.}

Ideals of minors of a matrix have been thoroughly studied over many decades. They play a significant role in the study of projective varieties. It had been a major classical problem to show that the ideal of $t \times t$ minors of a generic matrix is a prime and perfect ideal. The proof for a general value of $t$ was due to Eagon and Hochster from their important work in \cite{h-e}. In the first part of this paper, we generalize the notions of a matrix and its ideal of $2\times2$ minors to that of a {\it box-shaped-matrix} and its {\it ideal of $2\times2$ minors}. Our main theorem in this section is the following theorem.

\begin{thm}[Theorem \ref{1-generic-prime}]
If $\A$ is a box-shaped matrix of indeterminates, then $I_2(\A)$ is a prime ideal in $\goth{k}[\A]$ (here, $I_2(\A)$ is the ideal of $2\times2$ minors of $\A$).
\end{thm}

Coupled with previous work of Grone (\cite{grone}), we also show that the ideal of $2\times2$ minors of a {\it generic box-shaped matrix} is the defining ideal of a Segre embedding of the product of several projective space, namely $\mP(V_1) \times \ldots \times \mP(V_n) \hookrightarrow \mP(V_1 \otimes \ldots \otimes V_n)$. This geometric realization of the ideal of $2\times2$ minors of a generic box-shaped matrix enables us to study its perfection (Theorem \ref{1-generic-perfection}), its Hilbert function (Proposition \ref{1-hilbert-function}), and gives a Gr\"obner basis (Theorem \ref{1-generic-basis}). 

Box-shaped matrices not only describe the Segre embedding of the product of several projective spaces, but also provide a new tool for the study of projective embeddings of certain blowup surfaces. This study is carried out in the second part of this paper. To be more precise, let $\x = \{ P_1, \ldots$, $P_s\}$ be a set of $s$ distinct points in $\mP^2$, and let $\ix = \oplus_{t \ge \alpha} I_t \subseteq R = \goth{k}[w_1,w_2,w_3]$ be the homogeneous decomposition of the defining ideal of $\x$, and $\mP^2(\x)$ the blowup of $\mP^2$ centered at $\x$. The second part of this paper studies the problem of finding systems of defining equations for $\mP^2(\x)$ embedded in projective spaces by very ample divisors which correspond to the linear systems of plane curves going through the points in $\x$. This problem has also been considered by several authors in the last ten years, such as \cite{ge-gi}, \cite{ge-gi-h}, \cite{ge-gi-pi}, \cite{gi1}, \cite{gi2}, \cite{hol}, \cite{hol1} and \cite{hol2}.

A great deal of work has concentrated on an important special case, when $s={d+1 \choose 2}$ for some positive integer $d$ and the points in $\x$ are in generic position (cf. \cite{ge-gi}, \cite{gi1}, \cite{gi2}). In this case,
\[ \ix = I_d \oplus I_{d+1} \oplus I_{d+2} \oplus \ldots \]
is generated by $I_d$ (see \cite{g-m}). We also address this situation. 

It is well known that, in our situation, all the linear systems $I_t$ (for $t \ge d+1$) are very ample (cf. \cite{d-g}, \cite{ge-gi-pi}), so they all give embeddings of $\mP^2(\x)$ in projective spaces. If in addition, there are no $d$ points of $\x$ lying on a line, then the linear system $I_d$ is also very ample. Under this assumption, Gimigliano studied the embedding of $\mP^2(\x)$ given by the linear system $I_d$, which results in a White surface (\cite{gi1} and \cite{gi2}). White surfaces had also been studied in the classical literature (\cite{ro} and \cite{wh}). Gimigliano showed that the defining ideal of a White surface is generated by the $3 \times 3$ minors of a $3 \times d$ matrix of linear forms, and its defining ideal has the same Betti numbers as that of the ideal of $3 \times 3$ minors of a $3 \times d$ matrix of indeterminates (which was given by the Eagon-Northcott complex). The embedding of $\mP^2(\x)$ given by the linear system $I_{d+1}$ (which results in a Room surface) was then studied in detail by Geramita and Gimigliano (\cite{ge-gi}). Geramita and Gimigliano were able to determine the resolution of the ideals defining the Room surfaces. They also proved that the defining ideals of the embeddings of $\mP^2(\x)$ given by the linear systems $I_t$ are generated by quadrics, for all $t \ge d+1$, but they were unable to write down those generators when $t \ge d+2$. 

In \cite{chtv} and \cite{stv}, a method of finding a system of defining equations for a diagonal subalgebra from that of a bigraded algebra was given. This, together with results of \cite{m-u} (which gives the equations for the Rees algebra of the ideal of a set of ${d+1 \choose 2}$ points in generic position), makes it possible, in theory, for one to write a system of defining equations for the embeddings of $\mP^2(\x)$ given by the linear systems $I_t$ for all $t$. However, this method has its disadvantages as pointed out in \cite{ha-thesis}. In the second part of this paper, we generalize Geramita and Gimigliano's argument on the Room surfaces and give an explicit description of the defining ideals of the embeddings of $\mP^2(\x)$ given by the linear systems $I_t$, for all $t \ge d+1$. Our main result in this section is the following theorem.

\begin{thm}[Theorem \ref{2-main}] 
Suppose $t = d+n$ ~ ($n \ge 1$). Then the projective embedding of $\mP^2(\x)$ given by the linear system $I_t$ is generated by ${n+1 \choose 2}d$ linear forms and the $2 \times 2$ minors of a box-shaped matrix of linear forms. 
\end{thm}

Throughout this paper, $\goth{k}$ will be our ground field. For simplicity, we assume that $\goth{k}$ is algebraically closed and of characteristic 0 (though many of our results are true for any field $\goth{k}$). We also let $\mP^2 = \mP^2_{\goth{k}}$ be the projective plane over $\goth{k}$.

\section{Box-shaped matrices and their ideals of $2\times2$ minors.}

The techniques we use in this section are inspired by those of \cite{sh} in his study of ideals of $2\times2$ minors of a matrix.

\begin{center}
{\bf Box-shaped matrices}
\end{center}

Let $S$ be a commutative ring that contains a field $\goth{k}$. An $n$-dimensional array ($n \ge 2$)
\[ \A = (a_{i_1 \ldots i_n})_{1 \le i_j \le r_j, \ \forall j=1,\ldots,n}, \]
can be realized as the box
\[ \B = \{ (i_1,\ldots,i_n) | 1 \le i_j \le r_j, \ \forall j \}, \]
in which each integral point $(i_1,\ldots,i_n)$ is assigned the value $a_{i_1 \ldots i_n}$.

\begin{mydef} An $n$-dimensional array $\A$, with its box-shaped realization $\B$, is called an {\it $n$-dimensional box-shaped matrix} of size $r_1 \times \ldots \times r_n$.
\end{mydef}

We associate to each box-shaped matrix $\A$ of elements in $S$ a ring $S_{\A} = \goth{k}[\A]$, the subring of $S$ obtained by adjoining the elements of $\A$ to the field $\goth{k}$.

\begin{mydef} Suppose $\A$ is an $n$-dimensional box-shaped matrix of size $r_1 \times \ldots \times r_n$ of elements in $S$. For each $l = 1,2,\ldots,n$, we call
\[ a_{i_1 \ldots i_l \ldots i_n} a_{j_1 \ldots j_l \ldots j_n} - a_{i_1 \ldots i_{l-1} j_l i_{l+1} \ldots i_n} a_{j_1 \ldots j_{l-1} i_l j_{l+1} \ldots j_n} \in S_{\A}, \]
(where $(i_1,\ldots,i_n)$ and $(j_1,\ldots,j_n)$ are any two integral points in $\B$), a {\it $2\times2$ minor about the $l$-th coordinate} of $\A$. A {\it $2\times2$ minor} of $\A$ is a $2\times2$ minor about at least one of its coordinates. We let $I_2(\A)$ be the ideal of $S_{\A}$ generated by all the $2 \times 2$ minors of $\A$, and call it the {\it ideal of $2 \times 2$ minors} of the box-shaped matrix $\A$.
\end{mydef}

From now on, unless stated otherwise, we focus our attention to box-shaped matrices of indeterminates. Suppose $\A=(x_{i_1 \ldots i_n})_{(i_1,\ldots,i_n) \in \B}$ is an $n$-dimensional generic box-shaped matrix of size $r_1 \times \ldots \times r_n$ with its box-shaped realization $\B$. For each $l=1,\ldots,n$, let
\[ \A_l = (x_{i_1 \ldots i_n})_{(i_1,\ldots,i_n) \in \B \mbox{ and } i_l < r_l}, \]
and denote by $I_2(\A_l)$ its ideal of $2\times2$ minors (in the appropriate ring). For each $l=1,\ldots,n$, we also let
\[ \B_l = \{ (i_1,\ldots,i_n) \in \B | i_l = r_l \}, \]
and
\[ I_l = \ <I_2(\A_l), \{ x_{i_1 \ldots i_n} | (i_1,\ldots,i_n) \in \B_l \}>. \]

Throughout this paper, to any box-shaped matrix $\A$, we always associate box-shaped matrices $\A_l$, boxes $\B_l$ and all the ideals $I_l$ defined as above. The first crucial property of box-shaped matrices of indeterminates comes in the following lemma.

\begin{lemma} \label{1-intersection}
Suppose $\A=(x_{i_1 \ldots i_n})_{(i_1,\ldots,i_n) \in \B}$ is a box-shaped matrix of indeterminates in $S$. Then,
\renewcommand{\labelenumi}{(\alph{enumi})}
\begin{enumerate}
\item For any $l \not= s \in \{ 1,\ldots,n \}$, we have
\[ I_l \cap I_s = \ <I_2(\A), \{ x_{i_1 \ldots i_n} | (i_1,\ldots,i_n) \in \B_l \cap \B_s \}>. \]
\item For any distinct elements $l_1,l_2,\ldots,l_t$ of $\{ 1,2,\ldots,n \}$ ($2 \le t \le n$), we have
\[ \cap_{j=1}^t I_j = \ <I_2(\A), \{ x_{i_1 \ldots i_n} | (i_1,\ldots,i_n) \in \cap_{j=1}^t \B_j \}>. \]
\end{enumerate}
\end{lemma}

\begin{proof} {\bf (a)} For convenience, we denote by $LHS$ and $RHS$ the left hand side and the right hand side of the presented equality, respectively. It is clear that $RHS \subseteq LHS$. We need to show the opposite direction. Let $F \in LHS$. Since $F \in I_{l}$, we can write $F = F' + F''$, where
\[ F' \in I_2(\A_l), \mbox{ and } F'' = \sum_{(i_1,\ldots,i_n) \in \B_l} F_{i_1 \ldots i_n} x_{i_1 \ldots i_n}. \]
It suffices to show that $F'' \in RHS$. $F''$ certainly belongs to $I_s$. Now, for $(i_1,\ldots,i_n) \in \B_l$, we write $F_{i_1 \ldots i_n}$ in the form
\[ F_{i_1 \ldots i_n} = \sum_{(j_1,\ldots,j_n) \in \B_s} G_{i_1 \ldots i_n, j_1 \ldots j_n} x_{j_1 \ldots j_n} + G_{i_1 \ldots i_n}, \]
where $G_{i_1 \ldots i_n}$ is independent of the indeterminates $\{ x_{j_1 \ldots j_n} | (j_1,\ldots,j_n) \in \B_s \}$. Then $F'' = G + G'$, where
\[ G = \sum_{(i_1,\ldots,i_n) \in \B_l} G_{i_1 \ldots i_n} x_{i_1 \ldots i_n}, \]
and
\begin{eqnarray*}
G' & = & \sum_{(i_1,\ldots,i_n) \in \B_l, (j_1,\ldots,j_n) \in \B_s} G_{i_1 \ldots i_n, j_1 \ldots j_n} x_{i_1 \ldots i_n} x_{j_1 \ldots j_n} \\
& = & \sum_{\begin{minipage}{2.5cm}{\tiny $(i_1,\ldots,i_n)\in\B_l$, \\ $(j_1,\ldots,j_n)\in\B_s$} \end{minipage}} \bigg( G_{i_1 \ldots i_n, j_1 \ldots j_n} X_{i_1 \ldots i_n, j_1 \ldots j_n} + x_{i_1 \ldots i_{s-1} j_s i_{s+1} \ldots i_n} T_{i_1 \ldots i_n, j_1 \ldots j_n} \bigg),
\end{eqnarray*}
where
\[ X_{i_1 \ldots i_n, j_1 \ldots j_n} = x_{i_1 \ldots i_n} x_{j_1 \ldots j_n} - x_{i_1 \ldots i_{s-1} j_s i_{s+1} \ldots i_n} x_{j_1 \ldots j_{s-1} i_s j_{s+1} \ldots j_n}, \]
and $T_{i_1 \ldots i_n, j_1 \ldots j_n} \in S_{\A}$. Clearly, $X_{i_1 \ldots i_n, j_1 \ldots j_n}$ is a $2\times2$ minor about the $s$-th coordinate of $\A$, and since the sum is taken on $(i_1,\ldots,i_n) \in \B_l$ and $(j_1,\ldots,j_n) \in \B_s$, the point $(i_1,\ldots,i_{s-1},j_s,i_{s+1},\ldots,i_n)$ belongs to $\B_l \cap \B_s$. Thus, $G' \in RHS$.

It remains to show that $G \in RHS$. Again, we have $G \in I_s$, so we can write:
\[ G = H + \sum_{(j_1,\ldots,j_n) \in \B_s} H_{j_1 \ldots j_n} x_{j_1 \ldots j_n}, \]
where $H \in I_2(\A_s)$. We may also assume that $H$ and $H_{j_1 \ldots j_n}$, where $(j_1,\ldots,j_n) \in \B_l \cap \B_s$, are independent of the indeterminates
\[ \{ x_{j_1 \ldots j_n} | (j_1,\ldots,j_n) \in \B_s \backslash (\B_l \cap \B_s) \}. \]
Then
\[ G - H - \sum_{(j_1,\ldots,j_n) \in \B_l \cap \B_s} H_{j_1 \ldots j_n} x_{j_1 \ldots j_n} = \sum_{(j_1,\ldots,j_n) \in \B_s \backslash (\B_l \cap \B_s)} H_{j_1 \ldots j_n} x_{j_1 \ldots j_n}. \]
The left hand side of the above equality is independent of all the indeterminates
\[ \{ x_{j_1 \ldots j_n} | (j_1,\ldots,j_n) \in \B_s \backslash (\B_l \cap \B_s) \}. \]
Thus, both sides must be zero. This implies that
\[ G = H + \sum_{(j_1,\ldots,j_n) \in \B_l \cap \B_s} H_{j_1 \ldots j_n} x_{j_1 \ldots j_n} \subseteq RHS. \]

We have proved $LHS \subseteq RHS$. Thus, the given equality follows.

{\bf (b)} We will use induction on $t$. For $t = 2$ the equality is proved in part (a). Suppose $t > 2$, and the equality is true for $t-1$. We then have
\[ \cap_{j=1}^{t-1} I_j = \ <I_2(\A), \{ x_{i_1 \ldots i_n} | (i_1,\ldots,i_n) \in \cap_{j=1}^{t-1} \B_j \}>. \]

It remains to prove
\begin{eqnarray*}
\lefteqn{<I_2(\A), \{ x_{i_1 \ldots i_n} | (i_1,\ldots,i_n) \in \cap_{j=1}^t \B_j \}> \ =} \\
& = & <I_2(\A), \{ x_{i_1 \ldots i_n} | (i_1,\ldots,i_n) \in \cap_{j=1}^{t-1} \B_j \}> \cap \ I_t.
\end{eqnarray*}
We can proceed in the same lines of argument as that of the proof of part (a) to show that the equality above is indeed true. Hence, the presented equality is true for all $2 \le t \le n$.
\end{proof}

In particular, we obtain the following corollary.

\begin{cor} \label{1-new2}
$\cap_{l=1}^n I_l = \ <I_2(\A), x_{r_1 \ldots r_n}>$.
\end{cor}

\begin{center}
{\bf The prime-ideal theorem} 
\end{center}

From henceforth, we shall assume that our ring $S$ is a domain. The primeness of $I_2(\A)$ for a generic box-shaped matrix $\A$ comes as a consequence of a series of lemmas. Note that even though the following lemmas are generalizations to higher dimension of those given in \cite{sh}, most of the proofs require more arguments than what was given for their 2-dimensional statements.

\begin{lemma} \label{1-new3}
Suppose $F(\ldots,x_{i_1 \ldots i_n},\ldots)$ is an element of $S_{\A}=\goth{k}[\A]$. If, for some $x_{i_1 \ldots i_n}$ of $\A$, there exists a positive integer $\lambda$ such that $x_{i_1 \ldots i_n}^{\lambda} F \in I_2(\A)$, then, for any $x_{j_1 \ldots j_n}$ of $\A$ there exists a non-negative integer $\nu$ such that $x_{j_1 \ldots j_n}^{\nu} F \in I_2(\A)$.
\end{lemma}

\begin{proof}
Denote by $Z$ the multiplicatively closed subset of $S_{\A}$ consisting of all non-negative powers of $x_{j_1 \ldots j_n}$, and let $S_Z$ be the localization of $S_{\A}$ at the set $Z$. Let $\phi : S_{\A} \rightarrow S_Z$ be the ring homomorphism defined by $\phi(c) = c$ for all $c \in \goth{k}$, and
\[ \phi(x_{i_1 \ldots i_n}) = \sfrac{x_{i_1 j_2 \ldots j_n}x_{j_1 i_2 j_3 \ldots j_n} \ldots x_{j_1 \ldots j_{n-1} i_n}}{x_{j_1 \ldots j_n}^{n-1}} \mbox{ for all } x_{i_1 \ldots i_n} \in \A. \]

Obviously, $\phi$ is a well-defined map. It is easy to verify that $\phi(a) = 0$ for any $2\times2$ minors $a$ of $\A$. Thus, $\phi(I_2(\A)) = 0$. Moreover,
\[ x_{i_1 \ldots i_n}^{\lambda} F(\ldots,x_{i_1 \ldots i_n},\ldots) \in I_2(\A). \]
Therefore, in $S_Z$,
\[ \bigg( \sfrac{x_{i_1 j_2 \ldots j_n}x_{j_1 i_2 j_3 \ldots j_n} \ldots x_{j_1 \ldots j_{n-1} i_n}}{x_{j_1 \ldots j_n}^{n-1}} \bigg)^{\lambda} F\bigg(\ldots, \phi(x_{i_1 \ldots i_n}), \ldots\bigg) = 0. \]

Since $S_{\A}$ is a domain, so is $S_Z$. Hence,
\[ F\bigg( \ldots, \phi(x_{i_1 \ldots i_n}), \ldots \bigg) = 0 \mbox{ in } S_Z. \]

Now, using binomial expansions, we can write
\[ F(\ldots,x_{i_1 \ldots i_n},\ldots) = F\bigg( \ldots, \sfrac{x_{i_1 j_2 \ldots j_n}x_{j_1 i_2 j_3 \ldots j_n} \ldots x_{j_1 \ldots j_{n-1} i_n}}{x_{j_1 \ldots j_n}^{n-1}}, \ldots \bigg) + K, \]
where $K$ belongs to the ideal of $S_Z$ generated by elements of the form 
\[ x_{i_1 \ldots i_n} - \sfrac{x_{i_1 j_2 \ldots j_n}x_{j_1 i_2 j_3 \ldots j_n} \ldots x_{j_1 \ldots j_{n-1} i_n}}{x_{j_1 \ldots j_n}^{n-1}}. \]
The generators of this $S_Z$-ideal can be rewritten as 
\[ x_{i_1 \ldots i_n} x_{j_1 \ldots j_n}^{n-1} - x_{i_1 j_2 \ldots j_n}x_{j_1 i_2 j_3 \ldots j_n} \ldots x_{j_1 \ldots j_{n-1} i_n}. \]

We shall prove that $K$ belongs to the ideal of $S_Z$ generated by $I_2(\A)$, or equivalently, we prove that these generators, considered as elements of $S_{\A}$, belong to $I_2(\A)$. Indeed, using induction on $n$, modulo $I_2(\A)$, we have
\begin{eqnarray*}
K_n & = & x_{i_1 \ldots i_n} x_{j_1 \ldots j_n}^{n-1} - x_{i_1 j_2 \ldots j_n}x_{j_1 i_2 j_3 \ldots j_n} \ldots x_{j_1 \ldots j_{n-1} i_n} \\
& = & x_{i_1 j_2 \ldots j_n} x_{j_1 i_2 \ldots i_n} x_{j_1 \ldots j_n}^{n-2} + (x_{i_1 \ldots i_n} x_{j_1 \ldots j_n} - x_{i_1 j_2 \ldots j_n} x_{j_1 i_2 \ldots i_n}) x_{j_1 \ldots j_n}^{n-2} \\
&  & - x_{i_1 j_2 \ldots j_n}x_{j_1 i_2 j_3 \ldots j_n} \ldots x_{j_1 \ldots j_{n-1} i_n}  \\
& \equiv & x_{i_1 j_2 \ldots j_n} x_{j_1 i_2 \ldots i_n} x_{j_1 \ldots j_n}^{n-2} - x_{i_1 j_2 \ldots j_n}x_{j_1 i_2 j_3 \ldots j_n} \ldots x_{j_1 \ldots j_{n-1} i_n}  \\
& = & x_{i_1 j_2 \ldots j_n} K_{n-1},
\end{eqnarray*}
where
\[ K_{n-1} = x_{j_1 i_2 \ldots i_n} x_{j_1 \ldots j_n}^{n-2} - x_{j_1 i_2 j_3 \ldots j_n} \ldots x_{j_1 \ldots j_{n-1} i_n}. \]

Since every indeterminate appearing in the expression $K_{n-1}$ has $j_1$ in its first index, we can view $K_{n-1}$ as just the same expression as $K_n$ but given by the $(n-1)$-dimensional box-shaped matrix $\A' = (x_{i_1 \ldots i_n})_{i_1 = j_1}$. By induction hypothesis, $K_{n-1}$ then belongs to $I_2(\A') \subseteq I_2(\A)$. And hence, $K_n \in I_2(\A)$.

We have just proved that $K$ belongs to the ideal of $S_Z$ generated by $I_2(\A)$. Equivalently, $F(\ldots,x_{i_1 \ldots i_n},\ldots)$ belongs to the ideal of $S_Z$ generated by $I_2(\A)$. Therefore, there exists a $\nu$ such that
\[ x_{j_1 \ldots j_n}^{\nu} F(\ldots,x_{i_1 \ldots i_n},\ldots) \in I_2(\A) \mbox{ in } S_{\A}. \]
The lemma is proved.
\end{proof}

\begin{lemma} \label{1-new4}
Suppose $l \in \{ 1,2,\ldots,n \}$. Suppose also that $F \in S_{\A}=\goth{k}[\A]$ is a polynomial independent of the indeterminates $x_{i_1 \ldots i_n}$ for all $(i_1,\ldots,i_n) \in \B_l$ such that $I_2(\A_l) : F = I_2(\A_l)$. Then $I_l : F = I_l$.
\end{lemma}

\begin{proof} The proof follows in exactly the same line as that of \cite{sh}.
\end{proof}

\begin{lemma} \label{1-new5}
Let $F \in S_{\A}=\goth{k}[\A]$ and suppose that $x_{1 \ldots 1}^{\lambda}F \in I_2(\A)$ for some positive integer $\lambda$. Then $F \in I_2(\A)$. In other words, $I_2(\A) : x_{1 \ldots 1}^{\lambda} = I_2(\A)$.
\end{lemma}

\begin{proof}
We use induction on $n$. When $n=2$, the result follows from that of \cite{sh}. Suppose $n > 2$, $\A$ is an $n$-dimensional box-shaped matrix of indeterminates, and the lemma is true for any box-shaped matrices of lower dimension. We now use induction on $r_1+\ldots+r_n$. We may assume that $r_i \ge 2$ for all $i = 1,\ldots,n$ (since otherwise, $\A$ collapses to an $(n-1)$-dimensional box-shaped matrix, and the result follows from the induction hypothesis), and the lemma is true for any $n$-dimensional box-shaped matrix with smaller value of $r_1+\ldots+r_n$.

If $F$ is of degree zero, then $x_{1 \ldots 1}^{\lambda} F$ belongs to the ideal of $2 \times 2$ minors of a box-shaped matrix obtained from $\A$ by letting all the indeterminates $x_{i_1 \ldots i_n}$, for $(i_1,\ldots,i_n) \not= (1,\ldots,1)$, be zero. Yet, this ideal is zero, so $F = 0 \in I_2(\A)$. We may use induction again, assuming that the degree of $F$ is bigger than zero, and the lemma holds for polynomials whose degrees are smaller than that of $F$.

Now, $x_{1 \ldots 1}^{\lambda} F \in I_2(\A) \subseteq \cap_{j=1}^n I_j$ by Corollary \ref{1-new2}, so in particular, $x_{1 \ldots 1}^{\lambda} F \in I_j$ for all $j$. Moreover, by the induction hypothesis, we have $I_2(\A_j) : x_{1 \ldots 1}^{\lambda} = I_2(\A_j)$. Thus, by Lemma \ref{1-new4}, $I_j : x_{1 \ldots 1}^{\lambda} = I_j$. This implies that $F \in \cap_{j=1}^n I_j \ = \ <I_2(\A), x_{r_1 \ldots r_n}>$. Write $F = F_1 + x_{r_1 \ldots r_n} F_2$, where $F_1 \in I_2(\A)$. Since $I_2(\A)$ is homogeneous, we may assume that the degree of $F_2$ is smaller than that of $F$. We have $x_{1 \ldots 1}^{\lambda} F = x_{1 \ldots 1}^{\lambda} F_1 + x_{1 \ldots 1}^{\lambda} x_{r_1 \ldots r_n} F_2 \in I_2(\A)$. Thus, $x_{r_1 \ldots r_n} x_{1 \ldots 1}^{\lambda} F_2 \in I_2(\A)$. By Lemma \ref{1-new3}, there is a non-negative integer $\nu$ such that $x_{1 \ldots 1}^{\lambda + \nu} F_2 = x_{1 \ldots 1}^{\nu} x_{1 \ldots 1}^{\lambda} F_2 \in I_2(\A)$. By our induction hypothesis on the degree of $F$, we have $F_2 \in I_2(\A)$. Hence, $F \in I_2(\A)$ as required.
\end{proof}

The primeness of the ideal of $2 \times 2$ minors of a box-shaped matrix in the generic case is stated as follows.

\begin{thm} \label{1-generic-prime}
If $\A$ is a box-shaped matrix of indeterminates, then $I_2(\A)$ is a prime ideal in $\goth{k}[\A]$.
\end{thm}

\begin{proof}
Suppose that $F(\ldots,x_{i_1 \ldots i_n},\ldots) G(\ldots,x_{i_1 \ldots i_n},\ldots) \in I_2(\A)$, where $F, G \in S_{\A}=\goth{k}[\A]$. Let $Z$ be the multiplicatively closed subset of $S_{\A}$ consisting of all non-negative powers of $x_{1 \ldots 1}$, and let $S_Z$ be the localization of $S_{\A}$ at $Z$. Similar to what was done in Lemma \ref{1-new3}, we define a map
\[ \varphi : S_{\A} \rightarrow S_Z, \]
by sending $\goth{k}$ to $\goth{k}$, and sending $x_{i_1 \ldots i_n}$ to $\sfrac{x_{i_1 1 \ldots 1} x_{1 i_2 1 \ldots 1} \ldots x_{1 \ldots 1 i_n}}{x_{1 \ldots 1}^{n-1}}$ for all $x_{i_1 \ldots i_n} \in \A$. It is easy to verify that $\varphi(a) = 0$ for any $2\times2$ minors $a$ of $\A$. Thus, $\varphi(I_2(\A)) = 0$. Moreover, $F(\ldots,x_{i_1 \ldots i_n},\ldots) G(\ldots,x_{i_1 \ldots i_n},\ldots) \in I_2(\A)$. Hence, in $S_Z$,
\[ F(\ldots, \varphi(x_{i_1 \ldots i_n}),\ldots) \ G(\ldots, \varphi(x_{i_1 \ldots i_n}),\ldots) = 0. \]

Since $S_{\A}$ is a domain, so is $S_Z$. Thus, at least one of the two factors has to be zero. Suppose
\[ F(\ldots,\sfrac{x_{i_1 1 \ldots 1} x_{1 i_2 1 \ldots 1} \ldots x_{1 \ldots 1 i_n}}{x_{1 \ldots 1}^{n-1}},\ldots)  = 0. \]

Now, similar to what was done in Lemma \ref{1-new3}, we deduce that there exists a $\nu$ such that $x_{1 \ldots 1}^{\nu} F(\ldots,x_{i_1 \ldots i_n},\ldots) \in I_2(\A)$ in $S_{\A}$. Hence, by Lemma \ref{1-new5}, $F \in I_2(A)$, and this completes the proof.
\end{proof}

\begin{center}
{\bf Segre embedding, Cohen-Macaulayness and Kozsul property}
\end{center}

Suppose $V_1, V_2, \ldots, V_n$ are vector spaces of dimensions $r_1,r_2,\ldots,r_n$, respectively. Recall the following definition.

\begin{mydef}
A tensor $z \in V_1 \otimes \ldots \otimes V_n$ is referred to as {\it decomposable} if there exist $v_j \in V_j$ for all $j=1, \ldots, n$, such that $z = v_1 \otimes \ldots \otimes v_n$. 
\end{mydef}

Now, let $\{ e_{j1}, \ldots, e_{jr_j} \}$ be a basis for $V_j$ for all $j=1, \ldots, n$. Then a basis of $V_1 \otimes \ldots \otimes V_n$ is given by
\[ \{ \epsilon_{i_1 \ldots i_n} = e_{1 i_1} \otimes \ldots \otimes e_{n i_n} | 1 \le i_j \le r_j \ \forall j=1, \ldots, n \}. \]
A tensor $z \in V_1 \otimes \ldots \otimes V_n$ is represented by
\[ z = \sum_{i_1 \ldots i_n} y_{i_1 \ldots i_n} \epsilon_{i_1 \ldots i_n}, \]
and a vector $v_j \in V_j$ is given by
\[ v_j = \sum_{k=1}^{r_j} u_{j k} e_{j k}. \]
Thus, to have $z = v_1 \otimes \ldots \otimes v_n$, is the same as to have
\[ y_{i_1 \ldots i_n} = u_{1 i_1} \ldots u_{n i_n}, \mbox{ for all } i_1 \ldots i_n. \]
This is clearly the equations describing the image of the following Segre embedding:
\[ \mP(V_1) \times \ldots \times \mP(V_n) \hookrightarrow \mP(V_1 \otimes \ldots \otimes V_n). \]

Hence, a tensor $z \in V_1 \otimes \ldots \otimes V_n$ is decomposable if and only if its corresponding point in $\mP(V_1 \otimes \ldots \otimes V_n)$ is in the image of the above Segre embedding. 

The geometric realization of the ideal of $2\times2$ minors of a generic matrix $\A$ comes from the work of Grone (\cite{grone}), which we rephrase in the following proposition.

\begin{pro}[Grone, 1977]
Suppose $\A$ is a generic box-shaped matrix of size $r_1 \times \ldots \times r_n$, and $V_1, \ldots, V_n$ are vector spaces of dimension $r_1, \ldots, r_n$, respectively. Then $I_2(\A)$ gives a set of equations that describe the decomposable tensors in $V_1 \otimes \ldots \otimes V_n$.
\end{pro}

Since the Segre embedding of the product of several projective spaces is a closed immersion, Grone's result gives an immediate corollary, which demonstrates the geometric realization of $I_2(\A)$.

\begin{cor} \label{1-segre_embedding}
If $\A$ is an $n$ dimensional generic box-shaped matrix of size $r_1 \times \ldots \times r_n$, then $I_2(\A)$ gives the defining ideal of the Segre embedding
\[ \mP(V_1) \times \ldots \times \mP(V_n) \hookrightarrow \mP(V_1 \otimes \ldots \otimes V_n). \]
where $V_1, \ldots, V_n$ are vector spaces of dimensions $r_1, \ldots, r_n$, respectively.
\end{cor}

\begin{proof} The result follows from the fact that $I_2(\A)$ is a prime ideal.
\end{proof}

From this, we can calculate the Hilbert function of the ideal of $2\times2$ minors of a generic box-shaped matrix as follows.

\begin{pro} \label{1-hilbert-function}
The Hilbert function of $I_2(\A)$ is
\[ {\bf H}(I_2(\A), t) = {\prod_{i=1}^n r_i + t - 1 \choose t} - \prod_{i=1}^n {r_i + t-1 \choose t} \ \forall t \ge 0. \]
\end{pro}

\begin{proof} It is easy to see that all homogeneous polynomials of degree $t$ on $\mP^{\prod r_i -1}$ restricted to the image of $\mP^{r_1-1} \times \ldots \times \mP^{r_n-1}$ gives all multi-homogeneous polynomials of degree $(t, \ldots, t)$ in $\mP^{r_1-1} \times \ldots \times \mP^{r_n-1}$. Thus the Hilbert function of the homogeneous coordinate ring of the Segre embedding is $\prod_{i=1}^n {r_i + t-1 \choose t}$. The proposition now follows.
\end{proof}

{\bf Remark:} It is clear that any Segre embedding is Hilbertian, i.e. its Hilbert function and its Hilbert polynomial are the same.

The geometric realization of $I_2(\A)$ and Propostion \ref{1-hilbert-function} give us the perfection of $I_2(\A)$. The result is stated as follows.

\begin{thm} \label{1-generic-perfection}
If $\A$ is an $n$-dimensional generic box-shaped matrix of size $r_1 \times \ldots \times r_n$, then $I_2(\A)$ is a perfect ideal of grade $\prod_{i=1}^n r_i - \sum_{i=1}^n r_i + (n-1)$.
\end{thm}

\begin{proof}
We let $S_i = \goth{k}[y_{i,1},\ldots,y_{i,r_i}]$ be the homogeneous coordinate ring of $\mP^{r_i-1}$ for all $i$. Clearly, $S_i$ is Cohen-Macaulay for all $i$. By results of \cite[page 378]{vogel} and Proposition \ref{1-hilbert-function}, it follows by induction on $n$ that the Segre product $\underline{\otimes}_{i=1}^n S_i$ is a Cohen-Macaulay ring. Furthermore, this ring is exactly the coordinate ring of the Segre embedding $\mP^{r_1-1} \times \ldots \times \mP^{r_n-1} \hookrightarrow \mP^{\prod r_i -1}$. Thus, since $I_2(\A)$ is the defining ideal of this Segre embedding, i.e. $\underline{\otimes}_{i=1}^n S_i \simeq \goth{k}[\A]/I_2(\A)$, we have $I_2(\A)$ is a perfect ideal. The grade of $I_2(\A)$ comes from the codimension of the Segre embedding, which is exactly $\prod_{i=1}^n r_i - \sum_{i=1}^n r_i + (n-1)$. The theorem is proved.
\end{proof}

{\bf Remark:} The perfection of $I_2(\A)$ also comes from a more general result of Hochster (\cite[Theorem 1]{ho1}). 

We have an immediate corollary.

\begin{cor}
Suppose $V_1, \ldots, V_n$ are vector spaces of dimensions $r_1, \ldots, r_n$. Then, the homogeneous coordinate ring of the Segre embedding
\[ \mP(V_1) \times \ldots \times \mP(V_n) \hookrightarrow \mP(V_1 \otimes \ldots \otimes V_n) \]
is always Cohen-Macaulay.
\end{cor}

We now recall the following folklore result (cf. \cite{hu}).

\begin{lemma}[\cite{hu}, Lemma 6.3A] \label{1-generic}
Let $I$ be a proper ideal of $\z[x_1,\ldots,x_n]$ such that \linebreak $\z[x_1,\ldots,x_n]/I$ is $\z$-flat, and $I$ is perfect of grade $g$. Suppose $S$ is a Noetherian ring and $a_1,\ldots,a_n$ are elements of $S$. Let $I'$ be the ideal given by the image of $I$ under the ring homomorphism of $\z[x_1,\ldots,x_n] \rightarrow S$ sending $x_i$ to $a_i$. Then $\mbox{grade } I' \le g$, and if the equality is attained then $I'$ is a perfect ideal.
\end{lemma}

We note also that our calculations and arguments, so far, are independent of the field $\goth{k}$. In fact, the same calculations and arguments would apply if we have any commutative Noetherian ring with identity instead of $\goth{k}$. Thus, our results hold when we substitute $\goth{k}$ by $\z$, the ring of integers. This, together with Lemma \ref{1-generic}, gives rise to the following result for any box-shaped matrix $\A$.

\begin{thm} \label{1-general-perfection}
Suppose $\A$ is any $n$-dimensional box-shaped matrix of size $r_1 \times \ldots \times r_n$. Then,
\[ \mbox{grade } I_2(\A) \le \prod_{i=1}^n r_i - \sum_{i=1}^n r_i + (n-1), \]
and if the equality is attained then $I_2(\A)$ is a perfect ideal.
\end{thm}

\begin{proof}
The result follows from Lemma \ref{1-generic} and the fact that our Theorem \ref{1-generic-perfection} is still true if instead of $\goth{k}$ we have the ring $\z$.
\end{proof}

We return to the generic situation. Suppose again that
\[ \A=(x_{i_1 \ldots i_n})_{(i_1,\ldots,i_n) \in \B} \]
is a generic box-shaped matrix of size $r_1 \times \ldots \times r_n$. The following theorem gives a Gr\"obner basis for $I_2(\A)$.

\begin{thm} \label{1-generic-basis}
Under the degree reverse lexicographic monomial ordering on $S_{\A} = \goth{k}[\A]$, in which the variables $x_{i_1 \ldots i_n}$ are ordered by lexicographic ordering on their indices (assuming that $1 < 2 < \ldots < n$), the $2\times2$ minors of $\A$ form a Gr\"obner basis for $I_2(\A)$.
\end{thm}

\begin{proof} Let $\lex$ be the lexicographic ordering on $\nn^n$. We order the variables of $S_{\A}$ by 
\[ x_{i_1 \ldots i_n} \le x_{j_1 \ldots j_n} \Leftrightarrow (i_1, \ldots, i_n) \ \lex \ (j_1, \ldots, j_n), \]
and use degree reverse lexicographic ordering on the monomials of $S_{\A}$. We shall prove that under this monomial ordering, the $2\times2$ minors of $\A$ form a Gr\"obner basis for $I_2(\A)$. 

Let ${\mathcal G}$ be the collection of all $2\times2$ minors of $\A$. It suffices to show that the leading terms of ${\mathcal G}$ generate the leading term ideal of $I_2(\A)$. By contradiction, suppose $F \in I_2(\A)$, and $T$, the leading term of $F$, is not generated by the leading terms of ${\mathcal G}$. Clearly, from the nature of $I_2(\A)$, $T$ is a monomial with at least 2 different indeterminates. We consider a new partial ordering on the indeterminates of $S_{\A}$, defined by
\[ x_{i_1 \ldots i_n} \preceq x_{j_1 \ldots j_n} \Leftrightarrow i_l \le j_l \ \forall l=1, \ldots, n. \]

Suppose $x_{i_1 \ldots i_n}$ and $x_{j_1 \ldots j_n}$ are any two different indeterminates present in $T$. Without loss of generality, assume that $x_{i_1 \ldots i_n} < x_{j_1 \ldots j_n}$, i.e. there exists a positive integer $u$ such that $i_l = j_l$ for all $l=1,\ldots,u-1$, and $i_u < j_u$. It is easy to see that if $x_{i_1 \ldots i_n} \not\preceq x_{j_1 \ldots j_n}$ then there exists another integer $v > u$ such that $i_v > j_v$. In this case,   
\[ x_{i_1 \ldots i_{v-1} j_v i_{v+1} \ldots i_n} < x_{i_1 \ldots i_n}, x_{j_1 \ldots j_n} < x_{j_1 \ldots j_{v-1} i_v j_{v+1} \ldots j_n}. \]
Thus, $x_{i_1 \ldots i_n}x_{j_1 \ldots j_n}$ is the leading term of 
\[ x_{i_1 \ldots i_n}x_{j_1 \ldots j_n} - x_{i_1 \ldots i_{v-1} j_v i_{v+1} \ldots i_n}x_{j_1 \ldots j_{v-1} i_v j_{v+1} \ldots j_n} \in {\mathcal G}, \]
whence $T$ is generated by the leading terms of ${\mathcal G}$, a contradiction. Hence, these two indeterminates must be comparable, i.e. $x_{i_1 \ldots i_n} \preceq x_{j_1 \ldots j_n}$. This is true for any two different indeterminates of $T$. Therefore, $T$ can be rewritten as 
\[ T = x_{t_{11} \ldots t_{1n}}x_{t_{21} \ldots t_{2n}} \ldots x_{t_{p1} \ldots t_{pn}}, \]
for some positive integer $p \ge 2$, where
\[ x_{t_{11} \ldots t_{1n}} \preceq x_{t_{21} \ldots t_{2n}} \preceq \ldots \preceq x_{t_{p1} \ldots t_{pn}}. \]

Now, let $[y_{i,1}: \ldots: y_{i,r_i}]$ represent the homogeneous coordinates of $\mP^{r_i-1}$ for all $i=1, \ldots, n$. Since $I_2(\A)$ is the defining ideal of the Segre embedding 
\[ \mP^{r_1-1} \times \ldots \times \mP^{r_n-1} \hookrightarrow \mP^{\prod r_i - 1}, \]
$F$ vanishes when we substitute the indeterminate $x_{i_1 \ldots i_n}$ by $\prod_{l=1}^n y_{l, i_l}$ for all $(i_1, \ldots, i_n)$. It is also clear that after this substitution, $F$ becomes a polynomial on the variables $y_{i, j}$. This polynomial is zero for all values of the variables $y_{i, j}$, so it must be the zero polynomial (since the ground field $\goth{k}$ is infinite). This implies that there must be a term $T'$ of $F$ ($T' \not= T$) which cancels $T$ after the substitution. Suppose $x_{k_1 \ldots k_n}$ is an indeterminate present in $T'$. Since $T'$ cancels $T$ after the substitution, for each $l=1, \ldots, n$, $k_l \in \{ t_{1l}, \ldots, t_{pl} \}$. From the partial ordering on the indeterminates in $T$, it is now clear that $k_l \ge t_{1l}$ for all $l=1,\ldots,n$, whence $x_{t_{11} \ldots t_{1n}} \le x_{k_1 \ldots k_n}$. If $x_{t_{11} \ldots t_{1n}} < x_{k_1 \ldots k_n}$ for every indeterminate $x_{k_1 \ldots k_n}$ in $T'$, then $T < T'$, which is a contradiction since $T$ is the leading term of $F$. Otherwise, suppose $x_{t_{11} \ldots t_{1n}}$ is contained in $T'$, then by considering $T/x_{t_{11} \ldots t_{1n}}$ and $T'/x_{t_{11} \ldots t_{1n}}$, and continuing the process, we eventually would, again, get a contradiction. 

The theorem is proved.
\end{proof}

{\bf Remark:} From the proof above, it is easy to see that the $2\times2$ minors of $\A$ form a Gr\"obner basis for $I_2(\A)$ under any monomial ordering on $S_{\A}$ that satisfies the condition that if $g = x_{i_1 \ldots i_n}x_{j_1 \ldots j_n} - x_{p_1 \ldots p_n}x_{q_1 \ldots q_n}$ is an element of ${\mathcal G}$, where $x_{p_1 \ldots p_n} \preceq x_{q_1 \ldots q_n}$, then $x_{i_1 \ldots i_n}x_{j_1 \ldots j_n}$ is the leading term of $g$. Degree reverse lexicographic monomial ordering is merely one of those monomial orderings that satisfies this condition. We choose this ordering since it is pratical in most computational algebra packages, such as CoCoA and Macaulay2.

The theorem gives rise to an interesting corollary.

\begin{cor} \label{1-kozsul}
Suppose $V_1, \ldots, V_n$ are vector spaces of dimensions $r_1, \ldots, r_n$. Then, the homogeneous coordinate ring of the Segre embedding
\[ \mP(V_1) \times \ldots \times \mP(V_n) \hookrightarrow \mP(V_1 \otimes \ldots \otimes V_n) \]
is a Kozsul algebra.
\end{cor}

\begin{proof} This follows from the fact that all $2\times2$ minors of $\A$ are quadratic forms.
\end{proof}

\begin{center}
{\bf 3-dimensional box-shaped matrices} 
\end{center}

In the last part of this section, we briefly look at a particular class of box-shaped matrices, those of dimension 3. Besides the usual matrices, 3-dimensional box-shaped matrices are the easiest that can be visualized. To visualize all the $2\times2$ minors of a 3-dimensional box-shaped matrix, one only needs to take any two lines parallel to one of the axes, and looks at their intersection with any two planes parallel to the other two axes of our fixed system of coordinates. 3-dimensional box-shaped matrices not only describe the Segre embedding of the product of 3 projective spaces, but also give a tool in studying certain blowup surfaces, as it will be discussed in the next section. We first extend the notion of a box-shaped matrix of indeterminates to that of a {\it weak box-shaped matrix of indeterminates}.

\begin{mydef}
Suppose $\A=(a_{ijk})_{(i,j,k) \in {\bf B}}$ is a box-shaped matrix of forms in a ring $S$. For each integer $l$ let $\A_{(x,l)}$ be the matrix given by the collection $\{ a_{ijk} | (i,j,k) \in {\bf B}, i=l \}$. We call $\A_{(x,l)}$ an {\it x-section} of the box-shaped matrix $\A$. The y-sections and z-sections of $\A$ are defined similarly.
\end{mydef}

\begin{mydef}
A box-shaped matrix $\A = (x_{ijk})_{(i,j,k) \in {\bf B}}$, with box-shaped realization ${\bf B}$, of forms in a domain $S$ is called a {\it weak box-shaped matrix of indeterminates} if
\renewcommand{\labelenumi}{(\alph{enumi})}
\begin{enumerate}
\item All the entries in $\A$ are indeterminates of $S$, i.e. algebraically independent over $\goth{k}$.
\item $<I_2(\A),x_{r_1r_2r_3}> \ = \ \cap_{l=1}^3 I_l$ where the ideals $I_l$ are defined as that of a general $n$-dimensional box-shaped matrix.
\item There exists an integral point $(i,j,k) \in B$ such that when we set all indeterminates other than $x_{ijk}$ of $\goth{k}[\A]$ to zero, the ideal $I_2(\A)$ is the zero ideal.
\item The ideals of $2 \times 2$ minors of sections $\A_{x,i}$, $\A_{y,j}$ and $\A_{z,k}$ are prime ideals.
\end{enumerate}
\end{mydef}
\vspace{-6ex}
With this bigger class of box-shaped matrices, the primeness of their ideals of $2\times2$ minors still holds.

\begin{pro} \label{1-main}
$I_2(\A)$ is a prime ideal in $\goth{k}[\A]$ for any weak box-shaped matrix of indeterminates $\A$.
\end{pro}

\begin{proof} First, we can always re-arrange the indices such that $(i,j,k)$ becomes $(1,1,1)$. The proof now follows in the same lines as that of Theorem \ref{1-generic-prime}.
\end{proof}

\section{Projective embeddings of blowup surfaces.}

Let $\x \subseteq \mP^2$ be a set of $s={d+1 \choose 2}$ points ($d \in \z, d \ge 1$) that are in generic position. Let $\mP^2(\x)$ be the blowup of $\mP^2$ along the points of $\x$, and let $\ix = \oplus_{t \ge d}I_t \subseteq R = \goth{k}[w_1,w_2,w_3]$ be the defining ideal of $\x$. Let $\ov{\Lambda_t}$ be the surface obtained by embedding $\mP^2(\x)$ using the linear system $I_t$ ~($t = d+n, ~ n \ge 1$). In this section, we give an explicit description for a system of defining equations for $\ov{\Lambda_t}$ for any $t$. We start by a simple result, which could be of folklore.

\begin{lemma} \label{2-gen-2}
Suppose $S,R$ and $T$ are Noetherian commutative rings with identity, and $\phi : S \rightarrow R$ and $\psi : R \rightarrow T$ are surjective ring homomorphisms. Suppose also that $f_1,f_2,\ldots,f_n$ are generators for $\mbox{ker } \phi \subseteq S$ and $g_1,g_2,\ldots,g_m$ are generators for $\mbox{ker } \psi \subseteq R$. Let $p_j$ be a preimage of $g_j$ for all $j$, then $f_1,\ldots,f_n,p_1,\ldots,p_m$ give a set of generators for $\mbox{ker } (\psi \circ \phi) \subseteq S$.
\end{lemma}

\begin{proof} Clearly $f_i$'s and $p_j$'s are all in $\mbox{ker } (\psi \circ \phi)$. Moreover, if $x \in \mbox{ker } (\psi \circ \phi)$ then either $\phi(x) = 0$ or $\phi(x)$ is a linear combination of the $g_j$'s. The result is now trivial. 
\end{proof}

To proceed, it follows from \cite{g-m} that $\ix$ is minimally generated in degree $d$. By the Hilbert-Burch theorem, these generators are the $d \times d$ minors of a $d \times (d+1)$ matrix, say $\bL$, of linear forms :
\[ \bL = ( L_{ij} ), ~ L_{ij} \in R_1 \mbox{ for } i=1,2,\ldots,d \mbox{ and } j=1,2,\ldots,d+1. \]
In this notation,
\[ \ix = (F_1,\ldots,F_{d+1}), ~ F_i = (-1)^{i+1}\mbox{det}(\bL ~ \backslash ~ i^{\mbox{\footnotesize th}} \mbox{ column}). \]

For $\alpha=(\alpha_1,\alpha_2,\alpha_3)$, we write $w^{\alpha}$ for $w_1^{\alpha_1}w_2^{\alpha_2}w_3^{\alpha_3}$, and denote $|\alpha| = \alpha_1 + \alpha_2 + \alpha_3$. A system of generators of the vector space $I_t$ is given by ${n+2 \choose 2}(d+1)$ forms $w^{\alpha}F_j$ for $j=1,2,\ldots,d+1$ and $|\alpha| = n$.

Consider the rational map
\[ \varphi : \mP^2 \m \mP^{p}, \ p = {n+2 \choose 2} (d+1) - 1, \]
given by $\varphi(P) = [w^{\alpha}F_j]$ (we order the $\alpha$'s by lexicographic ordering with $w_1 > w_2 > w_3$). $\ov{\Lambda_t}$ embedded in $\mP^{p}$ is given by the closure of the image of $\varphi$.

Let $z_1 = w_1^n, z_2 = w_1^{n-1}w_2, \ldots, z_u=w_3^n$, where $u={n+2 \choose 2}$ (again, we arrange the terms in lexicographic order). We use homogeneous coordinates $[x_{ij}]_{1 \le i \le u, 1 \le j \le d+1}$ of $\mP^{p}$ such that
\begin{eqnarray}
\varphi([w_1:w_2:w_3]) & = & [x_{ij}], \mbox{ where } x_{ij} = z_iF_j. \label{2-coord-equation}
\end{eqnarray}

The vector space dimension of $I_t$ is $(n+1)d + {n+2 \choose 2}$, so there must be ${n+1 \choose 2}d$ dependence relations among the $w^{\alpha}F_j$'s. Those relations can be found as follows.

Let $\beta = (\beta_1,\beta_2,\beta_3)$ with $|\beta| = n-1$. For each $l=1,2,\ldots,d$, we have
\[ 0 = \mbox{ det } \left(
\begin{array}{c} w^{\beta}L_{l1} ~ w^{\beta}L_{l2} ~ \ldots ~ w^{\beta}L_{l d+1} \\ \bL \end{array}
 \right) = \sum_{j=1}^{d+1} L_{lj}w^{\beta}F_j. \]

Since $L_{lj}=\sum_{k=1}^3 \lambda_{ljk}w_k$, so by grouping similar terms, we get
\[ \sum_{|\alpha| = n, 1 \le j \le d+1} \mu_{l \alpha j}w^{\alpha}F_j = 0, ~ ~ \forall l=1,2,\ldots,d, \]
where 
\[ \mu_{l \alpha j} = \sum_{w^{\beta}w_k = w^{\alpha}} \lambda_{ljk} \mbox{ for each } l, \alpha \mbox{ and } j. \]
These are the dependence relations of the $w^{\alpha}F_j$'s. In terms of $z_i$'s, we can rewrite them as
\[ \sum_{i,j} \mu_{lij}z_iF_j = 0, ~ ~ \forall l=1,2,\ldots,d. \]
These give rise to the following equations:
\begin{eqnarray}
\sum_{1 \le i \le u, 1 \le j \le d+1} \mu_{lij}x_{i j} & = & 0, ~ ~ \forall l=1,2,\ldots,d. \label{2-linear-equations}
\end{eqnarray}

There are $d$ relations of the form (\ref{2-linear-equations}) for each $\beta$, and the number of such $\beta$'s is ${n+1 \choose 2}$. By abuse of notation, we denote the collection of all these ${n+1 \choose 2}d$ relations by (\ref{2-linear-equations}). The relations in (\ref{2-linear-equations}) would be independent relations if we can show that the ${n+1 \choose 2}d \times {n+2 \choose 2}(d+1)$ matrix {\bf E} of the coefficients $\mu_{lij}$ has maximal rank. Indeed, we shall use a similar argument to that given by Geramita and Gimigliano (\cite{ge-gi}).

\begin{lemma}
{\bf E} has maximal rank.
\end{lemma}

\begin{proof} We assume, without loss of generality, that none of the points of $\x$ is $P=[0:0:1]$, and that the first minor of $\bL$, $F_1$, does not vanish at $P$. Suppose
\[ \bL = A_1w_1+A_2w_2+A_3w_3, \]
where the $A_i$s have entries in the ground field. This means that $A_3$ has maximal rank $d$ (since $F_1(P) \not= 0$).

If we arrange the $\beta$s in lexicographic order with $w_1 > w_2 > w_3$ then {\bf E} would have the form
\[ {\bf E} = \left[
\begin{array}{cccccccc} A_1 & \ldots & A_2 & \ldots & A_3 & 0 & & \\
\ldots & A_1 & \ldots & A_2 & \ldots & A_3 & 0 & \\
\ldots & \ldots & \ldots & \ldots & \ldots & \ldots & \ldots & \ldots \\
\ldots & \ldots & \ldots & \ldots & \ldots & \ldots & \ldots & A_3 \end{array}
\right] \]
(The $A_3$ of the latter row is totally on the right of the $A_3$ of the former row). On each $A_3$, take $d$ columns that give a matrix $A_3'$ which has nonzero determinant. Putting them all together, we obtain a ${n+1 \choose 2}d \times {n+1 \choose 2}d$ matrix, which looks like the following:
\[ {\bf E'} = \left[
\begin{array}{cccc} A_3' & & & \\ & A_3' & & 0 \\ X & & \ddots & \\ & & & A_3' \end{array}
\right], \]
a lower-triangular matrix. Clearly, $\mbox{det} {\bf E'} = (\mbox{det} A_3')^{{n+1 \choose 2}} \not= 0$. Thus, the matrix {\bf E} has maximal rank.
\end{proof}

Obviously, on $\varphi(\mP^2 \backslash \x)$, the coordinates of the points satisfy the equations in (\ref{2-linear-equations}). These are the equations coming from the dependence relations of the $w^{\alpha}F_j$'s that we are looking for.

Consider the matrix
\[ M = \left[
\begin{array}{cccc} x_{11} & x_{12} & \ldots & x_{1 ~ d+1} \\
x_{21} & x_{22} & \ldots & x_{2 ~ d+1} \\
\ldots & \ldots & \ldots & \ldots \\
x_{u1} & x_{u2} & \ldots & x_{u ~ d+1} \end{array}
\right] \]

It is easy to see that the points of $\varphi(\mP^2 \backslash \x)$ satisfy all the $ 2 \times 2$ minors of $M$. Denote the collection of these equations by (**).

Moreover, on $\varphi(\mP^2 \backslash \x)$, each column of $M$ has the form :
\[ \left( \begin{array}{c} z_1F_j \\ z_2F_j \\ \ldots \\ z_uF_j \end{array} \right) \]
where $z_1 = w_1^n, \ldots, z_u = w_3^n$ for some point $[w_1:w_2:w_3] \in \mP^2 \backslash \x$. Clearly, the $z_i$'s satisfy the defining equations of the Veronese surfaces, which are known to be the $2 \times 2$ minors of certain Catalecticant matrices (see \cite{pu} for definition). Thus, on $\varphi(\mP^2 \backslash \x)$, the coordinates $x_{1j}, x_{2j},\ldots, x_{uj}$ satisfy the $2 \times 2$ minors of the Catalecticant matrix $\mbox{Cat}(1,n-1;3)$ of size $3 \times {n+1 \choose 2}$, for all $j=1,2,\ldots,d+1$. Denote the collection of these equations by (***).

From (\ref{2-coord-equation}), on $\varphi(\mP^2 \backslash \x)$, we have:
\[ x_{1j}/z_1 = x_{2j}/z_2 = \ldots = x_{uj}/z_u, ~ \mbox{ for all } j=1,2,\ldots,d+1. \]
This can be rewritten as a number of systems of equations, one for each $i=1,2,\ldots,u$
\[ \left\{
\begin{array}{ccc} x_{ij}/z_i & = & x_{1j}/z_1 \\ x_{ij}/z_i & = & x_{2j}/z_2 \\ & \ldots & \\
x_{ij}/z_i & = & x_{uj}/z_u \end{array} \right. \mbox{ for } j=1,2,\ldots,d+1. \]
Those relations give us, for each $i=1,2,\ldots,u$ :
\[ (S_i) \left\{ \begin{array}{ccc} x_{ij}z_1 - x_{1j}z_i & = & 0 \\ x_{ij}z_2-x_{2j}z_i & = & 0 \\
\ldots & & \\ x_{ij}z_u-x_{uj}z_i & = & 0 \end{array} \right. \mbox{ for } j=1,2,\ldots,d+1. \]

It is not hard to see that if the coordinates of $Q=[x_{ij}] \in \mP^{p}$ and $P=[z_i] \in \mP^{u-1}$ satisfy system $(S_i)$ for some $i$, where $z_i \not= 0$, then they satisfy systems $(S_i)$ for all $i$.

Before going further, we prove a similar proposition to that of \cite{ge-gi}.

\begin{pro} \label{uniqueness}
Let $Q=[x_{ij}]$ be a point on $\mP^p$, and suppose the coordinates of $Q$ satisfy equations in (**). Then there exists a unique $P = [z_1:\ldots:z_u] \in \mP^{u-1}$ such that the homogeneous coordinates of $P$ and $Q$ satisfy the systems $(S_i)$ for all $i$.
\end{pro}

\begin{proof}
Since the coordinates of $Q$ satisfy equations in (**), the matrix $M(Q)$ has rank 1, i.e. the rows of $M(Q)$ are all multiples of any nonzero row of $M(Q)$. Suppose the first row of $M(Q)$ is not identically zero (similar argument works for other rows). Then there exist $\nu_i$, for $i=2,\ldots,u$, such that
\[ x_{ij} = \nu_i x_{1j}, \mbox{ for all } j=1,2,\ldots,d+1. \]
We want $P \in \mP^{u-1}$ such that the coordinates of $P$ and $Q$ satisfy the systems $(S_i)$ for all $i$. We first consider such $P$ that the coordinates of $P$ and $Q$ satisfy $(S_1)$. This is the same as solving for $z_1,\ldots,z_u$ from $(S_1)$. The coefficients matrix becomes (projectively) a collection of :
\[ N_j = \left(
\begin{array}{ccccc} 0 & 0 & 0  & \ldots & 0 \\
-\nu_2 & 1 & 0 & \ldots & 0 \\
-\nu_3 & 0 & 1 & \ldots & 0 \\
\vdots & & &  \ddots & \\
-\nu_u & 0 & 0 & \ldots & 1 \end{array}
\right), \mbox{ for } j=1,2,\ldots,d+1. \]

Since $N_j$ is independent of $j$ and has rank exactly $u-1$, the system $(S_1)$ has exactly one projective solution, that gives a unique point $P \in \mP^{u-1}$. Moreover, this $P$ clearly has non-zero $z_1$ entry. Thus, the coordinates of $P$ and $Q$ satisfy $(S_i)$ for all $i$. Hence, $P$ exists and is unique.
\end{proof}

Let $\V$ be the algebraic set in $\mP^{p}$ defined by all the equations in (\ref{2-linear-equations}), (**) and (***). We have the following theorem.

\begin{thm} \label{2-asset}
$\V = \ov{\Lambda_t}$ as sets.
\end{thm}

\begin{proof} Clearly, $\varphi(\mP^2 \backslash \x) \subseteq \V$. Since $\V$ is closed, $\ov{\Lambda_t}$ is integral (so $\ov{\Lambda_t}$ is irreducible, and $\ov{\Lambda_t} = \overline{\varphi(\mP^2 \backslash \x)}$ ), we have
\[ \ov{\Lambda_t} \subseteq \V. \]
We only need to show that
\[ \V \subseteq \ov{\Lambda_t}. \]

Having Proposition \ref{uniqueness}, if we can show that for any points $P=[z_1:\ldots:z_u]$  and $Q=[x_{ij}]$ such that the coordinates of $Q$ satisfy the equations in (\ref{2-linear-equations}), (**) and (***), and coordinates of $P$ and $Q$ satisfy the systems $(S_i)$ for all $i$, $Q$ must be in $\ov{\Lambda_t}$, then we will have $\V \subseteq \ov{\Lambda_t}$, and so are done. Suppose $P$ and $Q$ are such points. We can always assume that $z_1 \not= 0$. Consider the system of equations given by all the equations in $(S_1)$ (if instead, $z_i \not= 0$, then we look at the system $(S_i)$). As a system of linear equations in indeterminates (note the way we have rearranged the indices)
\[ \{ x_{ij} | 1 \le j \le d+1, 1 \le i \le u \}, \]
the coefficients matrix is :
\[ A = \left(
\begin{array}{cccc} B & & & \\ & B & & \\ & & \ddots & \\ & & & B \end{array}
\right), \]
where
\[ B = \left(
\begin{array}{cccccc} 0 & 0 & 0 & \ldots & 0 \\
z_2 & -z_1 & 0 & \ldots & 0 \\
z_3 & 0 & -z_1 & \ldots & 0 \\
\vdots & & & \ddots & \\
z_u & 0 & 0 & \ldots & -z_1 \end{array}
\right). \]

Clearly, $B$ has rank $u-1$, and has a non-trivial solution $[z_1:\ldots:z_u]$. Therefore, the solution to $A$ must have the form :
\[ [x_{ij}] = [c_1z_1:c_1z_2:\ldots:c_1z_u:c_2z_1:\ldots:c_2z_u:\ldots:c_{d+1}z_1:\ldots:c_{d+1}z_u], \]
where $c_1,\ldots,c_{d+1}$ are constants not all zero, and the indeterminates are ordered by $1 \le j \le d+1$ and $1 \le i \le u$.

Now, since the coordinates of $Q$ also satisfy the equations in (***), which are the defining equations of Veronese surfaces, there exists a unique point $T=[w_1:w_2:w_3] \in \mP^2$ such that $z_1 = w_1^n, z_2 = w_1^{n-1}w_2, \ldots, z_u = w_3^n$. Thus,
\begin{eqnarray}
Q & = & [c_1w_1^n:\ldots:c_{d+1}w_3^n]. \label{2-break_1}
\end{eqnarray}

Lastly, the coordinates of $Q$ satisfy ${n+1 \choose 2}d$ equations in (\ref{2-linear-equations}), so
\[ \bL(T) \left[ \begin{array}{c} c_1 \\ c_2 \\ \ldots \\ c_{d+1} \end{array} \right] = \left[ \begin{array}{c} 0 \\ 0 \\ \ldots \\ 0 \end{array} \right]. \]

If $T \not\in \x$, then $\bL(T)$ has rank exactly $d$. Thus,
\begin{eqnarray}
\left[ \begin{array}{c} c_1 \\ c_2 \\ \ldots \\ c_{d+1} \end{array} \right]  & = & \rho \left[ \begin{array}{c} F_1(T) \\ F_2(T) \\ \ldots \\ F_{d+1}(T) \end{array} \right]. \label{2-break_2}
\end{eqnarray}
This implies that $Q \in \ov{\Lambda_t}$.

If $T \in \x$, then $\bL(T)$ has rank exactly $d-1$, so there is a 2-dimensional solution space, and these resulting $Q$s lie on a line of $\V$, which is one of the exceptional lines of $\ov{\Lambda_t}$.

Hence, we always have $Q \in \ov{\Lambda_t}$. We have proved that $\V = \ov{\Lambda_t}$ as sets.
\end{proof}

To continue our study, we let $S = \goth{k}[x_{ij}]$ be the homogeneous coordinate ring of $\mP^{p}$. Suppose ${\bf C}$ is the Catalecticant matrix $\mbox{Cat}(1,n-1;3)$ of indeterminates $\{ z_i \}_{0 \le i \le {n+2 \choose 2}}$. ${\bf C}$ is of size $3 \times {n+1 \choose 2}$. Consider the box ${\bf B}$ of size $(d+1) \times 3 \times {n+1 \choose 2}$. Let $\A$ be the box-shaped matrix obtained by assigning each integral point $(i,j,k)$ of ${\bf B}$ the indeterminate $x_{il}$ where $l$ is the integer such that $z_l$ is at the $(j,k)$-position in ${\bf C}$.

\begin{lemma} \label{2-weakbox}
$\A$ is a weak box-shaped matrix of indeterminates.
\end{lemma}

\begin{proof} Clearly, each x-section of $\A$ has its ideal of $2 \times 2$ minors as the defining ideal of a Veronese surface, so its ideal of $2 \times 2$ minors is a prime ideal. Also, each y-section and z-section of $\A$ is a matrix of indeterminates, whence whose ideal of $2 \times 2$ minors is also a prime ideal. Moreover, $x_{111}$ surely satisfies property (c) of $\A$ being a weak box-shaped matrix. It remains to show that
\[ <I_2(\A), x_{(d+1) 3 {n+1 \choose 2}}> \ = \ \cap_{l=1}^3 I_l. \]

For convenience, we let $r_1 = d+1, r_2 = 3, r_3 = {n+1 \choose 2}$, and consider $\A$ as a box-shaped matrix of size $r_1 \times r_2 \times r_3$. We shall first prove that
\[ I_2 \cap I_3 \ = \ <I_2(\A), \{ x_{i r_2 r_3} | i=1,\ldots,r_1 \}>. \]

The proof will go in the same line as that of part (a) of Lemma \ref{1-intersection}. Firstly, it is clear that $<I_2(\A), \{ x_{i r_2 r_3} | i=1,\ldots,r_1 \}> \ \subseteq I_2 \cap I_3$. It remains to show the other inclusion. Let $F \in I_2 \cap I_3$. Doing exactly as we did before, we end up with $F = F'+G'+G$, where $F',G' \in I_2(\A)$, and
\[ G = \sum_{ik} G_{ik}x_{ir_2k}, \]
where $G_{ik}$'s are independent of the variables $x_{ijr_3}$. Again, we have $G \in I_3$, so we can write
\[ G = H + \sum_{i,j} H_{ij}x_{ijr_3}, \]
where $H \in I_2(\A_3)$. We may assume that the $H_{ir_2}$'s are independent of all the variables $\{ x_{ijr_3} | j \not= r_2 \}$. By the nature of the $2\times2$ minors of $\A$ and the symmetry (in construction) of Catalecticant matrices, it can be seen that if a $2\times2$ minor of $\A$ has one indeterminate belonging to $\{ x_{ijr_3} | (i,j,r_3) \in \B \}$ then it must have at least two adjacent indeterminates belonging to $\{ x_{ijr_3} | (i,j,r_3) \in \B \}$. Thus, by re-grouping and rewriting, we can always assume that $H$ is also independent of the indeterminates $\{ x_{ijr_3} | (i,j,r_3) \in \B \}$. Now, clearly, $G = H + \sum_{i} H_{ir_2}x_{ir_2r_3} \in \ <I_2(\A), \{ x_{i r_2 r_3} | i=1,\ldots,r_1 \}>$. We have shown that $I_2 \cap I_3 \ = \ <I_2(\A), \{ x_{i r_2 r_3} | i=1,\ldots,r_1 \}>$.

It now follows in the same line of the proof of part (b) of Lemma \ref{1-intersection} that
\[ <I_2(\A), \{ x_{i r_2 r_3} | i=1,\ldots,r_1 \}> \ \cap \ I_3 = \ <I_2(\A), x_{r_1r_2r_3}>. \]

The lemma is proved.
\end{proof}

We obtain the main result of this section as follows.

\begin{thm} \label{2-main}
The subscheme $\ov{\Lambda_t}$ in $\mP^{p}$ is defined by ${n+1 \choose 2}d$ linear forms and the $2 \times 2$ minors of a box-shaped matrix of linear forms.
\end{thm}
\vspace{-3ex}
\begin{proof}
Let $S = \goth{k}[x_{ij}]$ be the homogeneous coordinate ring of $\mP^{p}$. Let $\A$ be the weak box-shaped matrix of indeterminates as above, and again, let $I_2(\A)$ be the ideal of $2 \times 2$ minors of $\A$ in $\goth{k}[\A]$. We also let $\bi$ be the ideal generated by $I_2(\A)$ and all the linear equations in (\ref{2-linear-equations}). Let $\cv$ be the subscheme of $\mP^{p}$ defined by $\bi$.

It is easy to see that $\bi$ contains all the equations in (\ref{2-linear-equations}), (**) and (***), so as sets, $\cv \subseteq \V$ (where $\V$ is the subvariety of $\mP^{p}$ defined by the equations in (\ref{2-linear-equations}), (**) and (***)).

Suppose now that $P=[\overline{w_1}:\overline{w_2}:\overline{w_3}] \in \mP^2 \backslash \x$ and $Q=[\overline{x_{ij}}] = \varphi(P)$. Let $\overline{z_1} = \overline{w_1}^n, \overline{z_2} = \overline{w_1}^{n-1}\overline{w_2}, \ldots, \overline{z_u} = \overline{w_3}^n$ then $\overline{x_{ij}}=\overline{z_i}F_j(P)$. Consider a $2 \times 2$ minors $a_{(K,L,M,N)}$ of $\A$ corresponding to the 4 points $K,L,M$ and $N$ in the box-shaped realization of $\A$. There are 3 possibilities for the tuple $(K,L,M,N)$.

{\bf Case 1.} $K=(i,j,k), L=(m,j,p), M=(m,n,p)$ and $N=(i,n,k)$ for some integers $i,j,k,m,n$ and $p$ (when the projections of $K,L,M,N$ on the $zx$-plane collapse to a line).

{\bf Case 2.} $K=(i,j,k), L=(m,j,k), M=(m,n,p)$ and $N=(i,n,p)$ for some integers $i,j,k,m,n$ and $p$ (when the projections of $K,L,M,N$ on the $yz$-plane collapse to a line).

{\bf Case 3.} $K=(i,j,k), L=(m,n,k), M=(m,n,p)$ and $N=(i,j,p)$ for some integers $i,j,k,m,n,$ and $p$ (when the projections of $K,L,M,N$ on the $xy$-plane collapse to a line).

By the construction of $\A$ and the fact that $[\overline{z_1}: \ldots: \overline{z_u}]$ is in the Veronese surface, i.e. it satisfies all the $2\times2$ minors of {\bf C}, it is easy to check that $Q$ satisfies the minors $a_{(K,L,M,N)}$. This is true for any $Q \in \varphi(\mP^2 \backslash \x)$ and any $2\times2$ minors $a_{(K,L,M,N)}$ of $\A$, so $\varphi(\mP^2 \backslash \x) \subseteq \cv$, whence $\V \subseteq \cv$.

We have shown that in all cases, $\V \subseteq \cv$. Hence, as sets, $\cv = \V = \ov{\Lambda_t}$.

Now, by Proposition \ref{1-main}, we know that $I_2(\A)$ is a prime ideal. Consider the following sequence of surjective ring homomorphisms:
\[ \goth{k}[x_{ij}] \stackrel{\phi}{\rightarrow} \goth{k}[w^{\alpha}t_j] \stackrel{\psi}{\rightarrow} \goth{k}[w^{\alpha}F_j], \]
defined in the obvious way; that is, both $\phi$ and $\psi$ send $\goth{k}$ to $\goth{k}$, and $\phi$ sends $x_{ij}$ to $w^{\alpha}t_j$ where $w^{\alpha}$ is labelled $z_i$, and $\psi$ sends $w^{\alpha}t_j$ to $w^{\alpha}F_j$.

We note that in proving equalities (\ref{2-break_1}) and (\ref{2-break_2}), we actually proved more. Firstly, the proof of (\ref{2-break_1}) and the fact that $I_2(\A)$ is a prime ideal imply that $I_2(\A)$ is the kernel of $\phi$. Secondly, the proof of (\ref{2-break_2}) shows that if we consider the equations in (\ref{2-linear-equations}) as polynomials over the $w^{\alpha}t_j$'s, then those polynomials are zero exactly when $t_j = F_j$ (since $t_j = F_j$ at all but a finite set of points $\x$). This implies that $\goth{k}[w^{\alpha}t_j]/\goth{a} \simeq \goth{k}[w^{\alpha}F_j]$, where $\goth{a}$ is the ideal generated by the images of the equations in (\ref{2-linear-equations}) through $\phi$. Thus, $\goth{a}$ is the kernel of $\psi$. Now, by Lemma \ref{2-gen-2}, we conclude that {\bf I} is the kernel of $\psi \circ \phi$. In other words, {\bf I} is the defining ideal of $\ov{\Lambda_t}$ embedded in $\mP^p$ (since the homogeneous coordinate ring of $\ov{\Lambda_t}$ embedded in $\mP^p$ is exactly $\goth{k}[w^{\alpha}F_j]$). The theorem is proved.
\end{proof}

{\bf Remark:} When $t =d+1$, our box-shaped matrix $\A$ collapses to be a normal matrix of size $3 \times (d+1)$, and the above result coincides with that obtained by Geramita and Gimigliano in \cite{ge-gi}.

\begin{small}
{\it Acknowledgement}. {\sf This paper is part of the author's PhD thesis. I would like to thank my research advisor A.V. Geramita for his inspiration and guidance.}
\end{small}

\end{document}